\documentclass[runningheads]{llncs}
\usepackage[T1]{fontenc}
\usepackage{graphicx}
\usepackage[mathlines]{lineno}
\usepackage{amsmath, amssymb, mathtools, hyperref}
\usepackage[noend]{algpseudocode}
\usepackage{algorithm}
\newcommand{\R}{\mathbb{R}}

\newcommand{\M}{\mathcal{M}}

\newcommand{\ga}{\gamma}
\newcommand{\Log}{\mathrm{Log}}
\newcommand{\Exp}{\mathrm{Exp}}

\newcommand{\U}{\mathcal{U}}
\newcommand{\grad}{\textnormal{grad}}
\newcommand{\vphi}{\varphi}
\DeclareMathOperator{\argmin}{arg\,min}
\usepackage{color}
\newcommand{\todo}[1]{\bgroup\color{red}#1\egroup}

\begin{document}
%
\title{Ridge Regression on Riemannian Manifolds for Time-Series 
Prediction\thanks{Extended version of conference paper presented at 
GSI 2025~\cite{nava_gsi2025}.}}
\titlerunning{Ridge Regression for Manifold-valued Time-Series}
%
\author{Esfandiar Nava-Yazdani\orcidID{0000-0003-4895-739X}}
\authorrunning{E. Nava-Yazdani}

%
\institute{Department of Visual and Data-Centric Computing, Zuse Institute Berlin\\ Berlin, Germany\\ 
\email{navayazdani@zib.de}\\
\url{https://www.zib.de/members/navayazdani}
}
\maketitle              
\begin{abstract}
We propose a natural intrinsic extension of ridge regression from Euclidean spaces to general Riemannian manifolds for time-series prediction. Our approach combines Riemannian least-squares fitting via Bézier curves, empirical covariance on manifolds, and Mahalanobis distance regularization. A key technical contribution is an explicit formula for the gradient of the objective function using adjoint differentials, enabling efficient numerical optimization via Riemannian gradient descent. We validate our framework through synthetic spherical experiments (achieving significant error reduction over unregularized regression) and hurricane forecasting.

\keywords{Hurricane \and Forecasting \and Mahalanobis \and Prediction \and Manifold-valued \and Tikhonov}\\
\textbf{MSC} 53C22 \and 53A25 \and 53B50 \and 62M10
\end{abstract}
\section{Introduction}
Manifold-valued time series arise naturally across diverse scientific domains: directions, bird migration, hurricane tracks on the sphere, rigid body rotations in SO(3), diffusion tensor evolution on the manifold of symmetric positive-definite matrices, and shape trajectories, for instance in Kendall's shape space. Traditional Euclidean methods for time-series analysis ignore the intrinsic geometric structure of these spaces, leading to artifacts such as non-unit quaternions when predicting rotations or non-positive-definite predictions in covariance modeling. Intrinsic methods that respect the geometry of data are essential for accurate and physically meaningful predictions.

We use intrinsically defined polynomials to model manifold-valued time-series. To this end, we utilize the extension of the ordinary de Casteljau algorithm and the resulting \emph{Bézier curves} to manifolds, presented in~\cite{Popiel2007BzierCA} and~\cite{nava2013casteljau}. The latter also proposes some basic properties and presents applications to motion of rigid body, space of positive definite matrices, construction of canal and developable surfaces, and curve design on implicit surfaces and polyhedra. Moreover, the mentioned work also extends the ordinary de Boor algorithm, and therewith, Bézier splines (also called composite Bézier curves) to manifolds. Note that, variational splines~\cite{varsplines} also offer flexible models, but are computationally more costly. The work~\cite{singh2013hierarchical} characterizes manifold-valued polynomials by vanishing higher-order covariant derivative. However, this approach also suffers from high computational cost (polynomials are determined by solving higher-order differential equations involving the curvature tensor). Due to evaluation via de Casteljau algorithm, Bézier polynomials offer high computational efficiency.

We employ a natural extension of the widespread \emph{least-squares regression} technique to estimate the best-fitting Bézier curves that model the trajectories. We recall that smoothing, dimensional reduction, and suppressing inconsistencies and noise are advantages of regression. For a comprehensive introduction and applications, we refer to~\cite{nava2019geomix,nava2020geo} and~\cite{nava2022HierCeoModel}, which employed geodesic regression (as a counter part to the ordinary linear one) and~\cite{nava2023hur,haniknava2024dc} for the general case of manifold-valued splines. In the applications, along with the regression, we compute the coefficient of determination to examine the model's adequacy. While geodesic regression extends linear regression to manifolds and variational splines~\cite{varsplines} offer flexible trajectory models, these methods lack regularization mechanisms for handling correlated or limited training data.

\emph{Ridge regression} is a popular adaptation of Tikhonov regularization to standard regression to estimate coefficients in the underlying model, especially when independent variables are highly correlated. This method has been used in many fields, including computational science and engineering and econometrics. It often provides prediction efficiency over simple linear regression with a significant number of model parameters and helps reduce overfitting that can result from model complexity. For an overview, we refer to~\cite{overviewridgereg2009,golub1999tikhonov} and~\cite{overviewridge2022}. 

The \emph{main contributions} of this work are as follows. We propose a natural intrinsic extension of the ordinary ridge regression from Euclidean spaces to arbitrary Riemannian manifolds, which incorporates the Mahalanobis distance and correlation matrix, and allows prediction of trajectories in the underlying manifold. We also derive an explicit formula for the gradient of the corresponding objective function, whose minimizer determines the model parameters. In addition, we validate our approach on synthetic spherical data, and hurricane tracks and intensities (wind speeds).

This work is organized as follows. The next section is devoted to mathematical preliminaries from Riemannian geometry, least-squares technique, the polynomial manifold-valued model, and our approach, the ridge regression. In Section~\ref{sec: apps}, we discuss applications including synthetic validation and hurricane forecasting, along with numerical results.

\section{Prediction for Trajectories in Riemannian Manifolds}
\label{sec: ridge}
This section details our regression model for representing trajectories using best-fitting Bézier curves, followed by our approach to prediction. For background on Riemannian geometry, we refer the reader to~\cite{doCarmo1992} and the brief introduction in Appendix~\ref{subsec: primer}. 

Let $M$ be a finite-dimensional Riemannian manifold with Riemannian distance $dist$, exponential map $\exp$ and local inverse, logarithmic map $\log$, respectively, and set $I:=[0,1]$. Moreover, we denote the tangent space of $M$ in $x$ by $T_xM$, and the number of control points of the best-fitting Bézier curves (polynomials) by $n$, set $\M:=M^n$ (power manifold), and $\U:=U^n$, where $U$ refers to a normal convex neighborhood (also called totally normal) in $M$. Moreover, we set $\ga (t;x,y):=\exp_x(t\log_xy)$ with $t\in I$ and $x,y\in U$.

\subsection{B\' ezier Polynomials and Regression}
To model manifold-valued time-series, we utilize the extension of the ordinary de Casteljau algorithm and the resulting Bézier curves to manifolds, presented in \cite{Popiel2007BzierCA} and \cite{nava2013casteljau}. The latter work not only introduces these geometric polynomials but also establishes fundamental properties and demonstrates applications across diverse domains: motion of rigid bodies in $\mathrm{SO}(3)$
, analysis in the space of positive definite matrices (crucial for diffusion tensor imaging and covariance modeling), construction of canal and developable surfaces in geometric design, and curve design on implicit surfaces and polyhedra. The versatility of this framework stems from its intrinsic formulation, which relies solely on the manifold's geodesic structure rather than embedding-dependent constructions.

The work \cite{singh2013hierarchical} characterizes manifold-valued polynomials by vanishing higher-order covariant derivatives, providing an alternative perspective rooted in differential geometry. However, this approach suffers from high computational cost, as polynomials are determined by solving higher-order differential equations involving the curvature tensor---a computationally intensive task, especially for manifolds with complex curvature structure. In contrast, due to evaluation via the de Casteljau algorithm, Bézier polynomials offer high computational efficiency through their recursive geometric construction, requiring only repeated geodesic interpolations without solving differential equations.

Moreover, we employ a natural extension of the widespread least-squares regression technique to estimate the best-fitting Bézier curves that model the trajectories. Regression techniques offer several crucial advantages: smoothing noisy observations, achieving dimensional reduction by representing complex trajectories through a small number of control points, and suppressing inconsistencies inherent in real-world measurements. These benefits are particularly valuable when dealing with heterogeneous data sources or measurements subject to varying levels of precision. For a comprehensive introduction and applications of regression on manifolds, we refer to \cite{nava2019geomix,nava2020geo} and \cite{nava2022HierCeoModel}, which employed geodesic regression (as a counterpart to ordinary linear regression in Euclidean spaces) for shape analysis and epidemiological applications, and~\cite{nava2023hur,haniknava2024dc} for the general case of manifold-valued splines applied to trajectory analysis and morphometrics. In our applications, along with the regression, we calculate the determination coefficient ($R^2$) to examine the suitability of the model and ensure that the polynomial representation captures a substantial part of the variability in the observed trajectories.

Let $b_1, \cdots, b_n$ be points in $U$. These define by the de Casteljau algorithm a polynomial $p(t;b)$ with control points $b_1, \cdots, b_n$. Here $b$ stands for $(b_1,\cdots, b_n)$ and $n-1$ is the degree of the polynomial $p(.;b)$. In the following, we briefly describe the extension of the de Casteljau algorithm and B\' ezier curves to Riemannian manifolds, and for a comprehensive introduction and some applications in data science, refer to~\cite{nava2023hur,haniknava2024dc}.

\textbf{Geometric Intuition.} In Euclidean spaces, Bézier curves are constructed through repeated linear interpolation between control points. The elegance of the de Casteljau algorithm lies in its purely geometric construction: to evaluate a Bézier curve at parameter $t$, we recursively interpolate between consecutive control points, with each level of recursion reducing the number of points by one until a single point remains. This geometric construction generalizes naturally to Riemannian manifolds by replacing linear interpolation with geodesic interpolation.

\subsubsection{De Casteljau Algorithm} 
\begin{figure}[t]
        \centering
        \includegraphics[width=.65\linewidth]{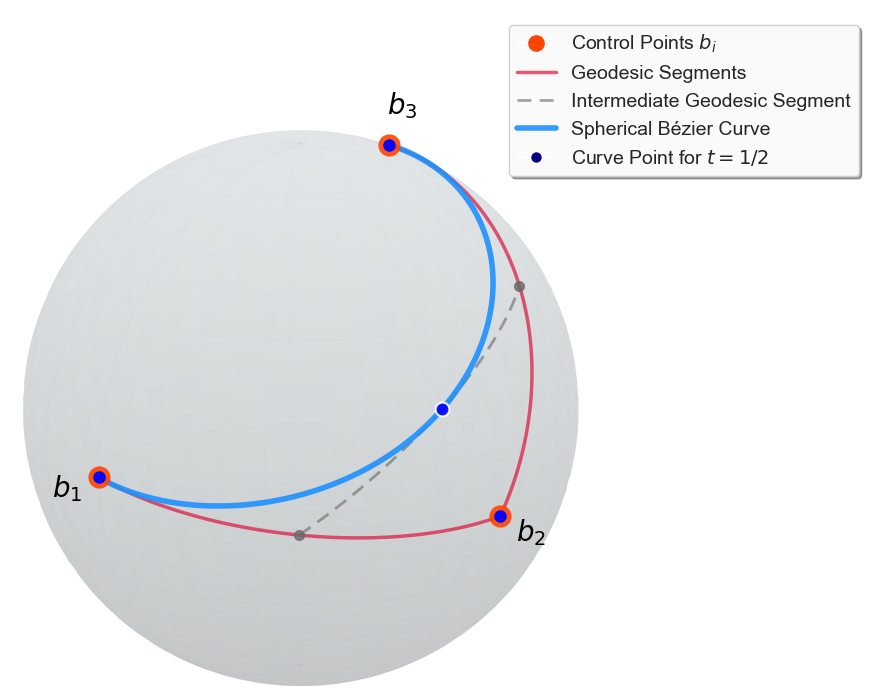}
        \caption{Construction of a quadratic Bézier curve on sphere via the de Casteljau algorithm. Control points $b_1, b_2, b_3$ (red) determine the curve (black). For parameter $t\in I$, we recursively interpolate: connect $b_1,b_2$ and $b_2,b_3$ with geodesics (blue), evaluate at $t$ to get intermediate points (green), then connect these with another geodesic and evaluate at t to obtain the curve point $p(t;b)$. This geometric construction generalizes naturally to higher degrees and arbitrary manifolds.}
        \label{fig: Bezier_curve_sphere}
        \end{figure}

For $t\in I$, we set
    \begin{align*}
       	&p_i^0(t):=b_{i+1}, \\
        &p_i^r(t):=\gamma(t; p_i^{r-1}(t), p_{i+1}^{r-1}(t)), \quad r=1,\dots,n-1, \quad  i=0,\dots,n-r.
    \end{align*}
We call $p(.;b) := p^{n-1}_0 : I \to U$ \textnormal{Bézier curve of degree $n-1$ with control points} $b_1,\dots,b_n$. Figure~\ref{fig: Bezier_curve_sphere} shows a visualization of the algorithm to construct a quadratic Bézier curve on the 2-sphere.

\textbf{Computational Efficiency.} Evaluation of Bézier polynomials is, due to the use of the de Casteljau algorithm, computationally very efficient. The algorithm requires $O(n^2)$ geodesic interpolations to evaluate the curve at a single parameter value, where each interpolation involves computing the exponential and logarithmic maps. For symmetric spaces, including spheres, hyperbolic spaces, and spaces of symmetric positive definite matrices, these maps have closed-form expressions, making the evaluation extremely fast. Even for general manifolds where geodesics must be computed numerically, the recursive structure of the algorithm remains stable and efficient. Evaluation of Bézier polynomials is, due to the use of the de Casteljau algorithm, computationally very efficient.

\subsubsection{Regression Model} Next, we use the widespread least-squares regression technique to estimate the best-fitting (not restricted to interpolatory) Bézier curves that model the trajectories in $M$.

Consider the scalars $t_1 < t_2 < \cdots < t_m$ and sample points $y_1,\cdots,y_m\in U$ of a trajectory $y$, i.e., $y_i=y(t_i)$. Our regression model aims to find a Bézier polynomial $p = p(.;b)$ that best fits the data $(t_i,y_i)$ in the least squares sense:
\[
\sum_{i=1}^mdist^2(y_i,p(t_i; b))\to\min.
\]
The minimizer is the best-fitting polynomial of degree $n-1$ for data $(t_i,y_i)$ with $i=1,\cdots,m$. Applying a linear bijection, we may assume that $t_1=0$ and $t_m=1$. Computationally, we employ the parametrization by the control points to determine the corresponding polynomial $p(.;b^\ast)$, where $b^\ast := \argmin H$ and
\[
H(b):=\sum_{i=1}^mdist^2(y_i,p(t_i;b)),\: b\in\U
\]
The choice $b_i=\ga(\frac{i-1}{n-1}; y_1,y_m)$ with $i=1,\cdots,n$ serves as a convenient initial guess. We recall that the $R^2$ value (coefficient of determination)
\begin{equation*}
R^2 = 1 - \frac{H_{\min}}{G_{\min}}
\end{equation*}
measures the goodness of fit (closer to 1 = better fit), where 
\begin{equation*}
G(x):=\sum_{i=1}^mdist^2(y_i,x),\: x\in U.    
\end{equation*}
\subsection{Ridge Regression}
Next, we propose our extended ridge regression. This approach, adapting a Tikhonov regularization to our parametric regression, allows prediction building on the polynomial model function. To this end, we utilize the notions of covariance and Mahalanobis distance, extended from Euclidean spaces to general manifolds (cf. the survey~\cite{pennec2006intro}). The regularizer has the geometric interpretation of Mahalanobis distance to the fitted prior data.

We denote the Riemannian metric of $\M$ by $g$, and its distance, exponential, and logarithmic maps by $\mathrm{d}$, $\Exp$ and $\Log$, respectively. Now, fix $\lambda \geq 0$, $\mu\in \U$, let $S$ be a positive definite selfadjoint endomorphism of $T_\mu \M$ and consider the function $F_{\lambda}: \U\to \R$ defined by
\begin{align*}
	F_{\lambda}(b) &=H(b)+\lambda g_\mu(S\Log_\mu b,\Log_\mu b)\\
    &=\underbrace{ \sum_{i=1}^mdist^2(y_i,p(t_i;b))}_{\text{data fidelity term}} + \lambda \underbrace{g_\mu(S\Log_\mu b,\Log_\mu b)}_{\text{regularization term}}.\end{align*}
Clearly, we may write the regularization term as \[g_\mu (S\Log_\mu x,\Log_\mu x)=g_\mu (W\Log_\mu x,W\Log_\mu x)\] where $W$ is a positive definite endomorphism with $WW^t=S$ (e.g. the unique positive definite symmetric root of $S$) on $T_\mu \M$. $\sqrt{\lambda}W$ is in the Euclidean case called Tikhonov matrix. We seek to minimize $F_{\lambda}$ with $\mu$ the expected value, $S$ the precision, i.e., the inverse of the covariance $\Sigma$ of $b$, and a suitable ridge parameter $\lambda$.
 
\subsubsection{Linearization and Gradient}

Let $\grad\, F_{\lambda}$ denote the gradient of $F_{\lambda}$. If $M$ is Euclidean, then $p(t;b)= Xb$ with a matrix $X$ determined by $(t_1,\cdots, t_m)$, and stacking $y_1,\cdots,y_m$ as a vector $y$, a straightforward computation shows that
\[
\grad_bF_{\lambda}=X^T(Xb-y) + \lambda S(b-\mu),
\]
implying the following well-known explicit expression for the minimizer of $F_{\lambda}$
\[
b = (X^TX+\lambda S )^{-1}(X^Ty+\lambda S\mu).
\]
Generally, due to absence of an explicit analytic solution, the regression task has to be solved numerically. To this end, we employ a Riemannian steepest descent solver~\cite{riemopt}. In this regard, the main challenge is to calculate the gradient of $F_{\lambda}$. In~\cite[Sec.\ 4.2]{BergmanGousenbourger2018} it has been shown that $\grad\, H$ can be computed as the adjoint of the sum of certain Jacobi fields (in closed form, if $M$ is a symmetric space). We still need to determine the gradient of the squared Mahalanobis distance. For this purpose, we have the following result.
\begin{proposition}[Gradient of squared Mahalanobis distance]
	Consider the function $f$ defined on $\U$ by
	\[
	f(x)=g_\mu (S\Log_\mu x,\Log_\mu x).
	\]
	Fix $x\in \U$. Let $\ga$ be the geodesic given by $\ga (t)=\Exp_\mu (t\Log_\mu x)$ with $t\in I$. Then
	\[
	\grad_x f = 2(d_x\Log_\mu)^\dagger S\Log_\mu x ,
	\]
    where the superscript $^\dagger$ represents the adjoint. Moreover, let $k$ be the function defined by $k(x)=\Exp_\mu (W\Log_\mu x)$ on $\U$, and denote the tangent map of $k$ evaluated at $x$ by $d_xk$. Then
	\[ \grad_x f= -2(d_xk)^\dagger\Log_{k(x)} \mu . \]
\end{proposition}
\begin{proof}
	Fix $v\in T_x\M$. Denoting the differential map of $f$ at $x$ by $d_xf$, we have the following.
	\begin{align*}
		d_xf v &= g_\mu (Sd_x\Log_\mu v ,\Log_\mu x)+g_\mu (S\Log_\mu x ,d_x\Log_\mu v )\\
		&=2g_\mu (S\Log_\mu x ,d_x\Log_\mu v )\\
	\end{align*}
    and arrive at
	\begin{align*}
		d_xf v &= 2g_x ((d_x\Log_\mu)^\dagger S\Log_\mu x ,v ).
	\end{align*}
	Hence
	\[
	\grad_x f = 2(d_x\Log_\mu)^\dagger S\Log_\mu x.
	\]
Now, consider the function $r$ defined by $r=d^2(\mu, .)$ on $U$. We have $f(x)=r(k(x))$. Thus,
	\begin{align*}
		d_xf v &= d_{k(x)}r (d_xkv)\\
		&= g_{k(x)}(\grad_{k(x)}r, d_xkv)\\
		&=g_x((d_xk)^\dagger(\grad_{k(x)}r), v)
	\end{align*}
 The well-known formula $\grad_yr=-2\Log_y\mu$ with $y=k(x)$ completes the proof. 
\end{proof}
Note that we may compute the gradient using automatic differentiation tools such as PyTorch or JAX. Moreover, if $M$ is embedded in some Euclidean space, then the gradient can also be obtained simply by projecting the Euclidean one to the tangent space at the point of evaluation. Finally, recall that for symmetric spaces Riemannian building blocks, and particularly the logarithmic map, are available as closed form expressions, accelerating the computations. For the computation of the adjoint map, we refer to Appendix~\ref{adjoint}.
\subsubsection{Mahalanobis Distance and Fitting Data}

Let $B$ denote a set of samples $b^1,\cdots, b^l$ in $\M$ with mean $\mu$. The matrix $\frac{1}{l-1}\sum_{j=1}^l\Log_\mu b^j(\Log_\mu b^j)^T$ (the superscript $T$ represents the transpose) in the exponential chart at $\mu$ defines an endomorphism $\Sigma$ on $T_\mu\M$, which is independent of the coordinates (cf.~\cite{pennec2006intro}). Thus, $\Sigma$ extends the notion of ordinary Euclidean (sample) covariance matrix to general Riemannian manifolds. Its inverse $S:=\Sigma ^{-1}$ is the precision operator. The squared Mahalanobis distance between $B$ and an arbitrary point $x\in M$ denoted by $d_{Mah}$ reads
\begin{align*}
	d_{Mah}^2(B,x)&=g_\mu (S\Log_\mu x,\Log_\mu x).
\end{align*}
Thus, the Mahalanobis distance accounts for data distribution. Now, consider a set $Y$ of discrete trajectories (not necessarily with equal sample sizes) in $M$ and let $B\subset\M$ be the set of corresponding control points of the best-fitting polynomials for the trajectories in $Y$, obtained by minimizing $H$. We may write 
\begin{align*}
	F_{\lambda} (b) &=H(b) + \lambda d_{Mah}^2(B,b).
\end{align*}
Thus, $F_{\lambda}$ averages the least-squares regression using $H$ given in the previous subsection and the squared Mahalanobis distance to $B$.

\subsubsection{Parameter Optimization and Forecasting}
The ridge parameter $\lambda$ controls the amount of regularization, which aims to achieve a balance between preventing overfitting and accurately predicting the training set. Nonetheless, an excessively high value of $\lambda$ may lead to an underfitted model that is unable to adequately capture significant patterns within the dataset. Therefore, it is important to determine an optimal value for $\lambda$.

Now, fix $\lambda$ and initial samples $y_1,\cdots, y_m$ of a discrete trajectory. We aim to predict $y_{m+1}$. Our prediction reads $\hat{y}_{m+1}:=p(t_{m+1};b^\ast)$, where $b^\ast$ denotes the minimizer of $F_{\lambda}$. A common core issue in forecasting is that, in general, for coarse sampling as well as for long-term forecasting ($t_{i}-t_{i-1}$ much larger than the preceding time step), the error $d(\hat{y}_i,y_i)$, due to noise and error propagation can increase excessively. To account for such effects and fine tuning, we slightly extend our ridge regression by adding a simple and computationally fast averaging step as follows. To this end, we consider an extended predictor $P_\alpha$
\begin{align*}
P_\alpha(t; b, x) = \exp_{x}(\alpha \log_{x}(p(t; b)))
\end{align*}
with $\alpha\in [0,1]$ and replace $\hat{y}_{m+1}$ by $P_{\alpha} (t_{m+1};b^*,y_m)$. Note that for $\alpha=1$, $\hat{y}_{m+1}$ remains unchanged. Iteration over $m$ yields the prediction for the whole trajectory.

Now, our task is to find optimal values for the parameter $\lambda$ and $\alpha$. To this end, one can use the well-known methods in the Euclidean case such as cross-validation. Thereby, denoting the validation dataset by $Y_{val}$, optimal parameter values $(\lambda^\ast, \alpha^\ast)$ can be gained iteratively via minimization of the residual 
\[
\sum_{y=(y_1, \cdots, y_m)\in Y_{val}}\sum_{i=1}^m d^2(\hat{y}_i,y_i).
\]
Obviously, the choice of the prior and its splitting in $Y_{val}$ and $Y$ depends on the application and, particularly, on properties such as seasonality and relevance of the historical order of the trajectories. Moreover, using a larger prior requires more computation time and memory. However, it does not necessarily result in better prediction accuracy on the test dataset. However, it is unlikely to achieve satisfactory prediction accuracy with smaller $Y$ and $Y_{val}$. Note that the ridge parameter $\lambda$ controls the tradeoff between complexity and fit by penalizing solutions that are too far from the mean $\mu$. For details on the complete forecasting pipeline, see Appendix~\ref{subsec: alg}.

\section{Numerical Experiments}
\label{sec: apps}
We demonstrate the effectiveness of our approach through two validations: (1) controlled synthetic experiments on spherical trajectories where ridge regularization reduces prediction error compared to unregularized regression, and (2) real-world hurricane forecasting where our method achieves competitive performance using only historical data.

In our experiments, $M$ was the 2-sphere embedded in $R^3$. We recall that the Riemannian exponential and logarithmic maps for $S^2$ read
\begin{align*}
	\exp_xv &=\cos (\vphi)\cdot x+\frac{\sin (\vphi)}{\vphi}\cdot v,\\
	\Log_xy &=\vphi\frac{y-\langle x,y\rangle x}{\| y-\langle x,y\rangle x \|},
\end{align*}
where $x,y\in S^2$ with $dist(x,y)<\pi$, $v\in T_x S^2$ with $\| v\| = \vphi$ and $\exp_x 0=x$. The computed $3n \times 3n$ covariance matrix in ambient coordinates is singular by construction, as all tangent vectors $\Log_\mu b$ lie in the $2n$-dimensional space $T_\mu(S^2)^n$. We employed Principal Component Analysis for dimensionality reduction, retaining only principal components with positive eigenvalues, effectively capturing the data's intrinsic dimensionality. For comparison, simpler regularized estimation methods, specifically Ad Hoc Diagonal Loading~\cite{diagload2005performance} and the Ledoit-Wolf technique~\cite{LEDOIT2004365}, were also applied. In all cases, these alternative regularization methods yielded results statistically equivalent to those obtained using the PCA-derived effective covariance matrix, confirming the robustness of the results against the specific choice of covariance estimator for this dataset.
\subsection{Synthetic Validation}
\begin{figure}[h]
        \centering
        \includegraphics[width=.5\linewidth]{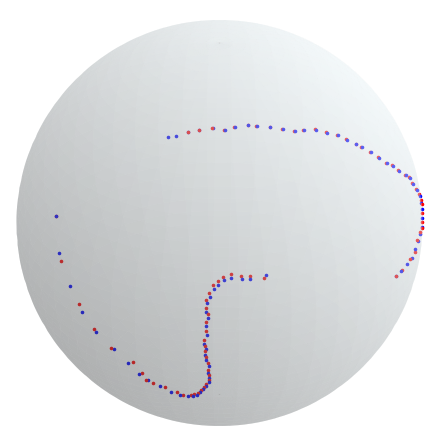}
        \caption{Two Representative forecasts from synthetic experiments, illustrating the prediction of the ridge regression (red) compared to the ground truth (blue). The close agreement demonstrates effective regularization.}
        \label{fig: synt}
    \end{figure}

A natural question is whether ridge regression on manifolds provides advantages compared to the ordinary least squares (OLS), $\lambda=0$, approach as known in the Euclidean case. Towards answering this question, we employed the following experiments with spherical data. We generated three template curves (mean trajectories): a geodesic, a cubic Bézier curve and the image of a perturbed sinusoidal curve in latitude-longitude coordinates, each consisting of $30--50$ points, and conducted 3 experiments ExpGeo, ExpCubic and ExpSin respectively. In each experiment, we generated 30 noisy trajectories by adding spatially-correlated random Gaussian noise in the tangent space at each point, with noise standard deviation $0.05$, followed by projection via the exponential map. This construction ensured between-trajectory correlation while maintaining realistic variability. Thus, in each case the resulting dataset comprised 30 trajectories, which were randomly split into 20 for training, 5 for validation and 5 for testing. We choose degree-5 Bézier polynomials yielding average $R^2$ value of $95\%$ to ensure high goodness of model-fit. 
For prediction, we set $\alpha =1$ and employed an iterative forecasting scheme where at each time step, we used the last $6$ observations to fit the model and predict the next position.
Using the optimal value $\lambda^*$ obtained by cross-validation yielded a minimum average MAE (mean absolute error) in all experiments smaller than the one from OLS ($\lambda = 0$), representing a $21.4\%$ average improvement. This demonstrates that ridge regularization on manifolds can provide substantial benefits when the data exhibit correlation structure, analogous to the Euclidean setting. Exemplary trajectories and their forecasts are depicted in Figure \ref{fig: synt}, and summary results are presented in Table \ref{tab: synt}.
        
\begin{table}[htbp]
\centering
\caption{Average MAE (radians) for synthetic spherical trajectory forecasting.}
\label{tab: synt}
\begin{tabular}{|l|c|c|c|}
\hline
Experiment & OLS ($\lambda=0$) & Proposed & Improvement \\
\hline
ExpGeo     & 0.035 & 0.028 & 20\% \\
\hline
ExpCubic   & 0.043 & 0.034 & 21\% \\
\hline
ExpSin     & 0.048 & 0.037 & 23\% \\
\hline
Average    & 0.042 & 0.033 & 21.4\% \\
\hline
\end{tabular}
\end{table}
\subsection{Application: Hurricane Forecasting}
Tropical cyclones, also known as hurricanes or typhoons, are among the most powerful natural phenomena with enormous environmental, economic, and human impact. The most common indicator of a hurricane's intensity is its maximum sustained wind speed, which is used to categorize the storm on the Saffir-Simpson hurricane wind scale. For example, wind speeds $\geq 137$ knots correspond to category 5. High track variability and out-most complexity of hurricanes has led to a large number of works to classify, rationalize and predict them. We remark that many approaches are not intrinsic and use linear approximations. We refer to the overview~\cite{OverviewMIHur2020}, the summary of recent progress~\cite{RecentProgress}, and~\cite{nhcGuid2022}.

In~\cite{nava2023hur}, we represented hurricane tracks using Bézier splines (compositions of Bézier polynomials), and computed average trajectories, performed tangent principal component analysis, and presented an intensity classification using support vector machines. In the following, we apply our proposed extended ridge regression approach, based on Bézier polynomial representation, to hurricane forecasting.

\subsubsection{Dataset}

\begin{figure}[htb]
	\centering
	\includegraphics[width=.55\textwidth]{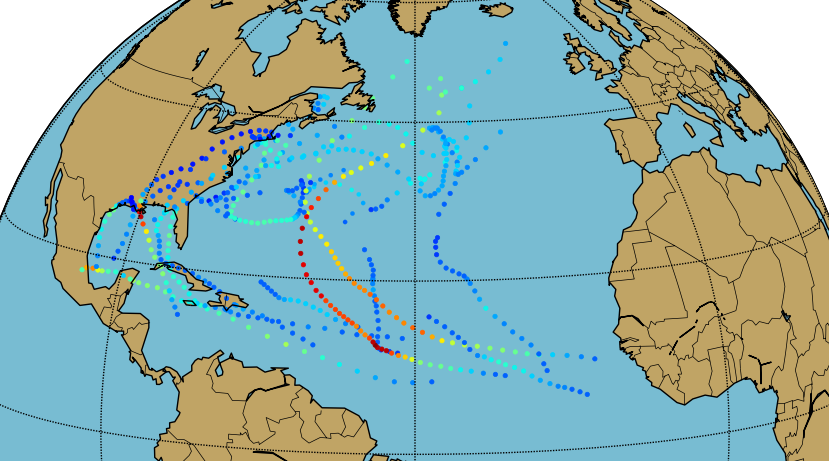}
	\hspace{20ex}
	\includegraphics[width=.55\textwidth,trim=0 0 10 0,clip]{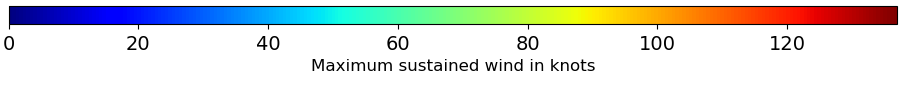}
	\caption{2021 Atlantic tracks and their intensities (maximum wind speeds).}
	\label{fig:hurricanes}
\end{figure}

We verify the effectiveness of the proposed framework by applying it to the Atlantic hurricane data from the publicly available \footnote{\url{https://www.nhc.noaa.gov/data/}} HURDAT 2 database provided by the U.S. National Oceanic and Atmospheric Administration. The data comprise measurements of latitude, longitude, and wind speed on a 6 hours base. The sample size (number of points that make up a track) varies from 13 to 96. Minimum of observed wind speeds is $15$ knots. Figure~\ref{fig:hurricanes} illustrates this dataset depicting the 2021 tracks together with their intensities.

\subsubsection{Experiments and Discussion}
\begin{figure}[htb]
	\centering
	\includegraphics[width=.55\textwidth,trim= 0 0 0 0,clip]{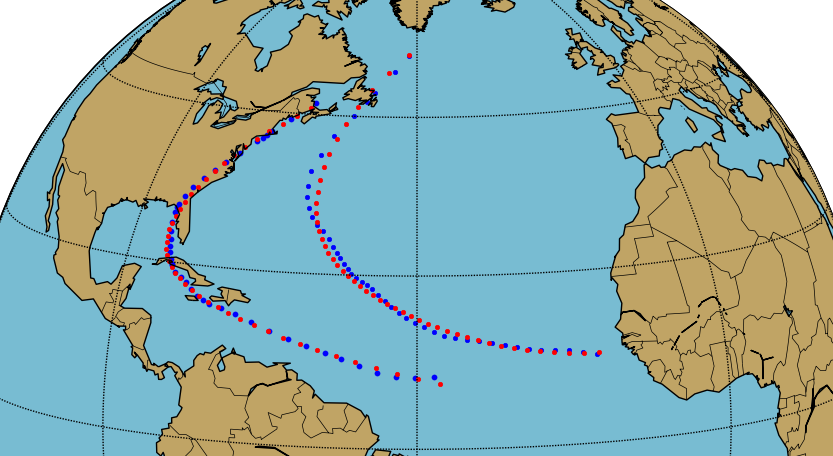}
	\caption{ Two hurricane tracks (blue) with their Bézier polynomial representations (red) using $n=6$ control points. High $R^2$ values indicate the polynomial model captures trajectory geometry effectively. These successful fits represent typical cases where our method performs well.}
	\label{fig:reghur}
\end{figure}

In our experiments, following the verification rules of the U.S. National Hurricane Center (NHC), we computed the average of the MAE (mean absolute error; using spherical distance for the tracks) of the 2021 Atlantic forecasts for comparison. To this end, we represented the tracks and intensities as discrete trajectories in a sphere of radius $3959$ miles (average earth radius) and $\R$, respectively. We initiated each forecast starting with the first sample and choose $n=6$, resulting in average $R^2$ values $> 0.95$ for both fitted tracks and intensities. The high $R^2$ values indicate that our polynomials fit well. Higher $n$-values although computationally more intensive, increased the coefficient of determination only incrementally by less then $1\%$. Figure~\ref{fig:reghur} shows two example tracks with their representation in our polynomial regression model.

We conducted two experiments, Exp1 and Exp2. Due to the deviations of annual trends and relevance of the chronological order, in Exp1, we used the 2020 ones (31 trajectories) as prior and the last 21 for validation to forecast all the 2021 ones (21 trajectories), and in Exp2, the first 16 trajectories from 2021 as prior with its last 5 trajectories to forecast the last 5 trajectories. It turned out that while the intensities are highly correlated, the covariance matrix of the tracks is almost singular (resolved as mentioned before). Higher values of $n$ slightly increased the forecast error of the tracks (overfitting due to small position correlations). For parameter optimization we only employed a simple grid search to decrease the computation time (on average a few seconds for each forecast). Variations in the settings with significantly larger or smaller validation sets slightly increased errors, and, as expected, validation sets chronologically closer to the test sets resulted in forecasts with slightly higher accuracies.

\begin{figure}[htb]
	\begin{center}
		\begin{minipage}{0.47\textwidth}
			\includegraphics[width=1\textwidth,trim=0 0 0 0,clip]{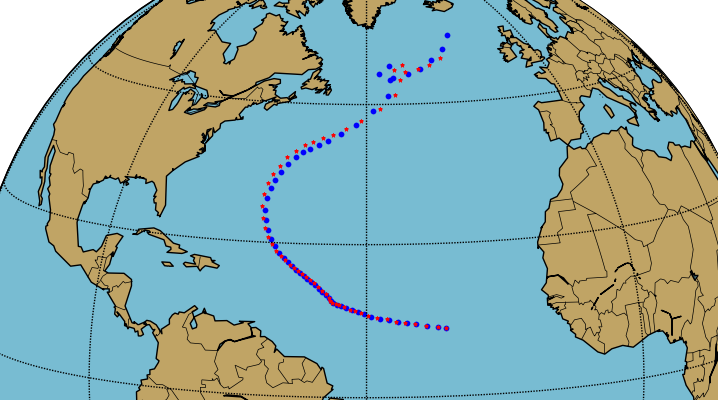}
		\end{minipage}
		\begin{minipage}{0.47\textwidth}
			\includegraphics[width=1\textwidth,trim=0 0 0 40,clip]{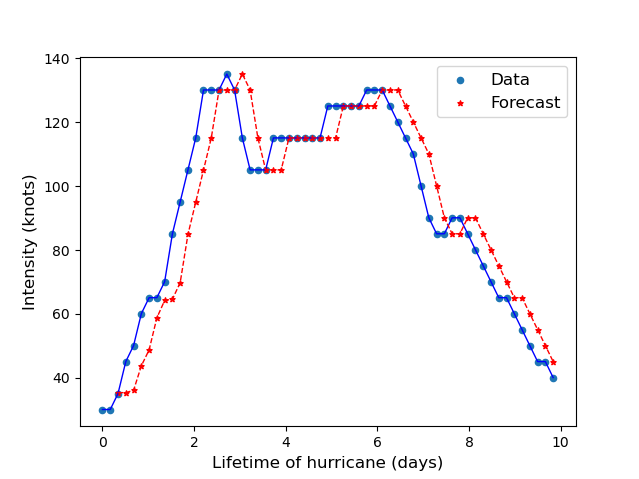}
		\end{minipage}
		\caption{Challenging forecast case: 12 h forecasts (red) via ridge regression for the category 4 and long-lasting 2021 anomalous hurricane Sam (blue). This hurricane has an anomalous highly complex dynamics, and very often and rapidly changes its intensity.}
		\label{fig:hurwind}
	\end{center}
\end{figure}

Obviously, one can also model the tracks as images of the Euclidean best-fitting polynomials in geographic coordinates (in the latitude-longitude plane) under the standard parametrization of the sphere. Our computations show that this increases errors ($5\%$ on average) demonstrating the advantage of our intrinsic geometry-aware approach. Figure~\ref{fig:hurwind} shows an example with a particularly large error, which is caused by the extreme complexity of dynamics and high variation, particularly near the loop. The increase in error due to the presence of a loop, in addition to meteorological anomalies, is an issue known in other approaches (cf. ~\cite{hurRNN2019pred}).
We did not expect our results to match the performance of the NHC's operational forecasts. Our geometry-aware approach is generic and relies purely on historical trajectory data and statistical-geometric regularization (via Mahalanobis distance), excluding the full complement of physical and atmospheric variables (such as vertical wind shear and sea surface temperature) that are explicitly incorporated into the NHC’s comprehensive suite of dynamical and statistical-dynamical models. Consequently, the NHC results serve as a state-of-the-art performance benchmark rather than a direct, like-for-like comparison. Table~\ref{tabhur} summarizes the average 12 h forecast errors of the standard NHC model OCD5 (cf.~\cite{nhc2021}) and our two primary experiments regarding the Atlantic basin for the 2021. Our geometric method achieves competitive performance with only a few training 
trajectories. Deep learning approaches like~\cite{hurRNN2019pred} require much more samples and lack 
interpretability. Our advantage is sample efficiency via geometric priors when training 
data is scarce (e.g., rare hurricane categories, emerging phenomena). Moreover, our framework currently treats position and intensity independently. A joint representation in $S^\times \R$ would capture position-intensity correlations (e.g., rapid intensification over warm water).

\begin{table}[htbp]
\centering
\caption{Average 12-hour forecast MAE for 2021 Atlantic hurricanes.}
\label{tabhur}
\begin{tabular}{|l|l|c|c|}
\hline
Experiment & Method & Intensity (kt) & Track (mi) \\
\hline
Exp1 & NHC       & 6.4 & 48.2 \\
\cline{2-4}
     & Proposed  & 6.9 & 90.7 \\
\hline
Exp2 & NHC       & 7.5 & 67.4 \\
\cline{2-4}
     & Proposed  & 6.5 & 97.2 \\
\hline
\end{tabular}
\end{table}
For the code implementing our approach, which in particular includes Riemannian optimization for the computation of geodesic paths and B\'ezier polynomials, we used the publicly available Python package \href{https://github.com/morphomatics}{morphomatics} v4.0 (cf.~\cite{Morphomatics}).

\section{Conclusions}\label{sec5}
In this work, we presented an intrinsic natural extension of polynomial ridge regression, respectively, Tikhonov regularization, from Euclidean spaces to arbitrary Riemannian manifolds. This allows for the prediction of manifold-value time series using the Mahalanobis distance to the training data as a prior. We also provided formulae for computing the gradient of the loss function in the approach. Furthermore, due to the use of the de Casteljau algorithm for temporal evaluation of our polynomial model, our approach is computationally fast. We also presented empirical clear error reduction on synthetical spherical data compared to the ordinary least-squared approach. In addition, we provided a discussion of our numerical experiments on the application of the approach to hurricane forecasting and intensity.

Beyond meteorology, our framework applies to: shape trajectories in medical imaging (organ growth, disease progression), rigid body motion in robotics (orientation tracking), covariance evolution in finance (portfolio dynamics), and morphometrics in biology (evolutionary shape changes). By providing explicit gradients and respecting intrinsic geometry, our method advances computational tools for geometric data science, offering sample-efficient alternatives to data-hungry machine learning approaches crucial when training data is scarce.

There are several exciting tasks for future work. First, we plan to consider further data-driven modifications to improve our results on hurricanes, especially with respect to long-term forecasts. Particularly, a joint representation in $S^2\times\R$ to incorporate position-intensity correlation is natural. We also plan to consider further applications and apply our approach to morphology and shape analysis. As another example application, and to exploit the effect of negative curvature, we plan to study the Hadamard manifold of positive definite symmetric matrices, which is of both application and theoretical importance.
\section{Appendix}
\subsection{Background: Riemannian Geometry}
\label{subsec: primer}
A {Riemannian manifold $M$ is a differentiable manifold endowed with a Riemannian metric $\langle \cdot, \cdot \rangle$, which is a smoothly varying inner product $\langle \cdot, \cdot \rangle_p$ on the tangent space $T_p M$ for every $p \in M$. This metric induces the Riemannian distance function $dist(p, q)$, defined as the infimum of lengths of curves between $p$ and $q$. If $f: M \to N$ is a smooth map between two manifolds $M$ and $N$, we denote the derivative of $f$ at $p$ in direction $v \in T_pM$ by $d_pf(v)$. It defines a linear mapping $d_pf: T_pM \to T_{f(p)}N$. The Levi-Civita connection $\nabla$ applied to vector fields $X$ and $Y$, $\nabla_XY $ is the directional derivative of $Y$ along $X$. The manifold's generalization of a straight line, is called geodesic (characterized by zero acceleration). A curve $\ga:I\to M$ is geodesic iff $\nabla_{\gamma'}\gamma' = 0$. For optimization, we rely on tools defined locally within a normal convex neighborhood $U \subset M$, where geodesics are unique and length-minimizing.\\

\textbf{Core Geometric Maps}
\begin{itemize}
    \item Exponential Map ($\exp_p$): Maps a tangent vector $v \in T_p M$ to the endpoint of the unique geodesic starting at $p$ with initial velocity $v$, i.e., $\exp_p(v) := q$.
    \item Logarithmic Map ($\log_p$): The inverse of $\exp_p$, which returns the initial velocity vector $v = \log_p(q) \in T_p M$ of the geodesic connecting $p$ to $q$.
\end{itemize}

\textbf{Calculus and Statistics}\\
The metric structure generalizes necessary calculus and statistical tools for optimization:
\begin{itemize}
    \item Gradient: For a smooth function $f: M \to \mathbb{R}$, the gradient at $p$ is defined implicitly by $\d f_p(v) = \langle \textnormal{grad}_p f, v \rangle_p$ for all $v \in T_pM$. This generalization retains the property that the gradient points in the direction of steepest ascent~\cite{riemopt}.
    \item Fréchet Mean ($\overline{q}$): The sample mean for data $q_1,\dots,q_m \in U$ is generalized as the point that minimizes the sum of squared distances (the Riemannian variance):
    \begin{equation*} \label{def:frechet_mean}
        \overline{q} := \argmin_{p \in U}\sum_{i=1}^m dist(p, q_i)^2.
    \end{equation*}
\end{itemize}
\subsection{Computation of Adjoint Map}
\label{adjoint}
Let $x,y\in M$, $v\in T_xM$ and $w\in T_yM$ and $F:T_xM\to T_y M$. We recall how to compute $F^\dagger$ (the adjoint of $F$).\\

\textbf{Using embedding:}\\
Let $A$ denote the matrix representation of $F$ w.r.t. an ONB of the ambient Euclidean space with inner product $\langle\cdot,\cdot\rangle $. Then
\begin{align*}
    g_y(Fv,w) & = g_x(v, F^\dagger w)\\
    \langle Fv, w\rangle & = \langle v, F^\dagger w\rangle\\
    (Av)^Tw &= v^T (A^\dagger w)\\
    v^TA^Tw &= v^TA^\dagger w
\end{align*}
Thus, $A^\dagger = A^T$.\\

\textbf{Using local coordinates in a chart:}\\
Let $G_x$ and $G_y$ denote the (Gram) matrix representations of $g$ in local charts at $x$ and $y$ respectively. Then
\begin{align*}
    g_y(Fv,w) & = g_x(v, F^\dagger w)\\
    (Av)^TG_yw &= v^T G_x (A^\dagger w)\\
    v^TA^TG_yw &= v^TG_x A^\dagger w
\end{align*}
Hence $G_x^{-1}A^TG_y w = A^\dagger w$ implying $A^\dagger = G_x^{-1}A^TG_y$.

\subsection{Algorithmic Forecasting Pipeline}
\label{subsec: alg}
The following summarizes the algorithm for prediction via proposed extended ridge regression.
The main objective is to find $b^*$ that minimizes the penalized least-squares functional $F_\lambda$.
Since an explicit solution is generally not available, the minimization is performed iteratively using a Riemannian optimization solver (e.g., gradient descent). We gain the training control points $B$ from training trajectories by fitting polynomials via least-square regression given by $H$. We compute optimal values for the ridge parameter $\lambda$ and the (fine tuning) averaging parameter $\alpha$ using cross-validation and minimizing the residual on a search grid in $[0,\infty[ \times [0,1]$.  
\begin{algorithm}[ht]
\caption{Forecasting via Extended Intrinsic Ridge Regression on Manifolds}
\label{alg:intrinsic_ridge_regression0}
\begin{algorithmic}[1]
\Statex \textbf{Preprocessing}
\Statex \textbf{1. Get training and validation sets:} Split prior trajectories in training and validation trajectories
\Statex \textbf{2. Training:} Fit Bézier polynomials to training trajectories $\rightarrow$ control points $B$
\Statex \textbf{3. Compute statistics:} Calculate Fréchet mean $\mu$ and precision $S$ from $B$
\Statex \textbf{4. Optimal $\lambda$ and $\alpha$:}
\Statex $\quad$ Require search grid $\Gamma=\{(0,0), \cdots, (\lambda_{max},1) \}\subset [0,\infty [ \times [0,1]$
\Statex $\quad$ Use cross-validation to find optimal values $\lambda^*$ and $\alpha^*$ by minimizing the residual on $\Gamma$

\Statex \textbf{Processing}
\Statex For a new history $(y_1,\cdots, y_m)$ predict $y_{m+1}$
\Statex \textbf{1. Optimal b via gradient descent}:  $b^* = \underset{b\in N}{\arg\min} F_{\lambda^* }(b)$
\Statex $\quad$ Require max iterations $K$, learning rate $\eta$, tolerance $\epsilon$
\Statex $\quad$ Initialize control points $b^{(0)}$ (e.g., set $b^{(0)} = \mu$)
\Statex$\quad$ for $k=0$ to $K-1$
    \Statex $\quad$ \text{Compute the gradient:}
    \Statex $\quad$ $G^{(k)} = \grad_b F_\lambda(b^{(k)})$
    
     if $\|G^{(k)}\| < \epsilon$
        \Statex$\quad$ break
    \Statex $\quad$ \text{Update the control points using the Riemannian exponential map:}
    \Statex $\quad b^{(k+1)} = \Exp_{b^{(k)}} \left( -\eta G^{(k)} \right)$
\Statex $\quad$ return $b^{*} = b^{(K)}$

\Statex \textbf{2. Result:} $\hat{y}_{m+1}=P_{\alpha^{*}}(t_{m+1}; b^{*}, y_m)$
\end{algorithmic}
\end{algorithm}

\section*{Acknowledgments}
The author is supported through the DFG individual funding with 
project ID 499571814. 

This paper significantly extends the conference version~\cite{nava_gsi2025} presented at the Geometric Science of Information (GSI) conference, 2025. The main additions include: (1) comprehensive synthetic validation experiments demonstrating significant error reduction, (2) enhanced covariance regularization analysis, (3) detailed algorithmic description, and (4) extended 
theoretical discussion and implementation details.
%
%
%
\bibliographystyle{splncs04}
\bibliography{lit}

\end{document}